\documentclass[11pt]{article}
\usepackage{rotating,multirow,amsmath,amssymb,fullpage}

\usepackage{natbib}
\bibpunct[, ]{(}{)}{,}{a}{}{,}%

\newtheorem{theorem}{Theorem}
\newtheorem{lemma}{Lemma}
\newtheorem{proposition}{Proposition}
\newenvironment{proof}{\paragraph{Proof:}}{\hfill$\Box$}

\newcommand{\mbfx}{\mathbf{x}}
\newcommand{\mbfy}{\mathbf{y}}
\newcommand{\mbfz}{\mathbf{z}}

\begin{document}
	\title{Solving Zero-sum Games using Best Response Oracles\\
		with Applications to Search Games}

	\author{%
		Lisa Hellerstein%
		\thanks{Department of Computer Science and Engineering, NYU Tandon School of Engineering, 6 Metrotech Center, Brooklyn, NY 11201, USA} 
		\and
		Thomas Lidbetter%
		\thanks{Department of Management Science and Information Systems, Rutgers Business School, Newark, NJ, USA}
		\and
		Daniel Pirutinsky%
		\footnotemark[2]
	}

	\maketitle

\begin{abstract} We present efficient algorithms for computing optimal or approximately optimal strategies in a zero-sum game for which Player I has $n$ pure strategies and Player II has an arbitrary number of pure strategies. We assume that for any given mixed strategy of Player I, a best response or ``approximate'' best response of Player II can be found by an oracle in time polynomial in $n$. We then show how our algorithms may be applied to several search games with applications to security and counter-terrorism. We evaluate our main algorithm experimentally on a prototypical search game.  Our results show it performs well compared to an existing, well-known algorithm for solving zero-sum games that can also be used to solve search games, given a best response oracle.
\end{abstract}

\section{Introduction}
\label{sec:intro}
Consider a zero-sum game with positive integer payoffs in which Player I (the maximizer) has $n$ pure strategies and Player II (the minimizer) has an arbitrary number of pure strategies. The payoff function is denoted by $C$, so that if Players I and II play pure strategies $i$ and $j$, respectively, then the payoff is given by $C(i,j)$. We also allow the arguments of $C$ to be probability vectors expressing mixed strategies, in which case $C$ denotes expected payoffs. Whenever we refer to games in this paper we assume they are of this form. Optimal mixed strategies and the value of the game can be found solving a linear program, but if Player II's strategy set is very large compared to $n$ then the computational time required to do this may not be polynomial in $n$.

Examples of this nature frequently arise in the study of finite search games for an immobile Hider. In such games, which are played between a Hider and a Searcher, the Hider typically has a strategy set corresponding to $n$ locations in which to hide and the Searcher's strategy set corresponds to a permutation of those locations (or a subset thereof), corresponding to the order in which she searches the locations. The order of search may be restricted due to a network structure of the search space or other restrictions on the mode of search. Notwithstanding this, the size of the Searcher's strategy set is usually exponential in $n$. The payoff of the game is generally either some cost incurred by the Searcher in searching for the Hider, which she wishes to minimize, or the probability of detection, which she wishes to maximize. 

Search games are motivated by wide-ranging problems in national security, counter-terrorism, search-and-rescue, biology and others. Summaries of the literature on search games for an immobile Hider can be found in \cite{AG03,Gal2011} and \cite{Hohzaki}.

In many such games, it is possible to exploit the structure of the game in order to solve the game (that is, find optimal strategies and the value) or find bounds on the value. An example of this is the classic search game for an immobile Hider on a network studied in \cite{Gal79} and \cite{Gal2000}. More recent examples can be found in \cite{Alpern16,AL15,AL13a,AL13b,ADL16,BK13,BK15,FKR16} and \cite{LS16}. Many of these games are infinite, in the sense that one or both of the players has a strategy set of infinite cardinality, but in this paper we restrict ourselves to finite games, which lend themselves better to an algorithmic approach.

Exploiting the structure of the game may not be possible in some games, or it may limit what can be achieved. However, in some games it may be possible to efficiently solve the one-sided problem of finding an optimal best response of Player II for any given mixed strategy of Player I: that is, finding a pure strategy for Player II which has minimal expected payoff against the given mixed strategy of Player I. We will refer to this problem as the {\em best response problem}. In this paper we study games for which there exists an oracle that can efficiently solve the best response problem. More generally, we consider games such that there is an oracle that, for any mixed strategy of Player I, can find a pure strategy for Player II that ensures the payoff is no more than $\alpha$ times the best response payoff. We call such an oracle an {\em $\alpha$-approximate best response oracle}. So if $\alpha=1$, then an $\alpha$-approximate best response oracle finds best responses; in this case we simply call it a {\em best response oracle}.

We show that for games with an $\alpha$-approximate best response oracle, there are efficient algorithms for finding $\alpha$-optimal strategies for both players. By this we mean a mixed strategy for Player I that guarantees a payoff of at least $V^*/\alpha$ and a mixed strategy for Player II that guarantees a payoff of at most $\alpha V^*$, where $V^*$ is the value of the game. 

Some of our analysis relies on previous work, as explained in Subsection~\ref{sec:overview}, and our main contribution is the wide range of applications to the field of search games, detailed in Section~\ref{sec:app}. The examples we give include known search games with known solutions, known search games with unknown solutions and new search games, all of which can be efficiently solved, either exactly or approximately, using our algorithms. We hope that by giving these techniques more prominence, they might be used by others for the purpose of solving search games, or indeed games in other fields. 

We demonstrate the efficacy of our approach experimentally by applying our main algorithm for computing optimal strategies (presented in Section 4) to the
game described in Subsection \ref{sec:boxes}. We show it performs favorably relative to a comparable algorithm of~\cite{FS99,FS96}.

\subsection{Example: Searching in Boxes}
\label{sec:boxes}
Before describing our algorithms in detail, we illustrate the type of game we will consider with an example of a search game that is easy to state, and that has the characteristics of the type of games to which our results can be applied. Although its solution is already known, it is a fundamental game and provides a good introduction to the more complex games we will discuss in Section~\ref{sec:app}.

An item (for example an explosive device) is hidden in one of $n$ boxes with costs $c_1,\ldots,c_n$, so that the Hider's strategy set is the set of boxes, $[n]=\{1,\ldots,n\}$. A pure strategy for the Searcher is simply an ordering of the boxes, so the Searcher's strategy set is the set of all permutations $\pi:[n] \rightarrow [n]$. In other words, $\pi(j)$ is the $j$th location to be inspected. Clearly the Hider's strategy set has size $n$ and the Searcher's strategy set has size $n!$. For given strategies $j$ of the Hider and $\pi$ of the Searcher, the payoff $C(j,\pi)$ of the game, which we call the {\em search cost} is the sum of the costs of all the boxes in the ordering up to and including the box that contains the hidden object. That is,
\[
C(j,\pi) = \sum_{i \le \pi^{-1}(j)} c_{\pi(i)}.
\]
We call this game BOX. The solution of BOX is already known in closed form. The Hider's optimal strategy is to choose box $j$ with probability proportional to $c_j$ as shown in \cite{AL13b}, where it was also proved that this is the unique optimal strategy. The Searcher does not have a unique optimal strategy, but three different optimal strategies can be found in \cite{AL13b,Lidbetter13} and (implicitly) in \cite{CDHW}. The value of the game is also given in closed form (in two different ways) in  \cite{AL13b} and \cite{Lidbetter13}. It is somewhat surprising, given the simplicity with which the game can be described, that prior to these publications the game had not been solved, to the best of our knowledge.

We now show that one could alternatively find optimal strategies in this game by applying our algorithms. To do this, we must show that there is a best response oracle. This is the problem of finding a permutation $\pi$ that minimizes $C(\mbfx,\pi)$ for a given Hider mixed strategy $\mbfx$. It is a classic search problem solved by Blackwell, reported in \cite{Matula}. In fact, Blackwell solved a more general problem, which we describe later in Subsection~\ref{sec:infinite}. For this version of the problem, the solution is to search the boxes in non-increasing order of the index $x_j/c_j$. 

We can arrive at the solution in another way, taken from the scheduling literature. Consider the problem of scheduling $n$ jobs with processing times $t_j$ and weights $w_j$, that correspond to the relative importance of the jobs. The object is to determine the schedule that orders the jobs in such a way as to minimize the weighted sum of the completion times. This problem is usually denoted $1||\sum w_j C_j$ in the scheduling literature, and has the well known solution, given by \cite{Smith}, of scheduling the jobs in non-decreasing order of the ratio $t_j/w_j$. Of course, viewing the costs $c_j$ in the best response problem for BOX as processing times, and viewing the probabilities $x_j$ as weights, these two problems and their solutions are equivalent.

We could also consider the more general box searching game introduced and solved by \cite{Lidbetter13} of searching for $k \ge 1$ objects hidden in $n$ boxes. However, in this case our approach cannot be used to solve the game, since the solution of the best response problem is not known. In fact, the computational complexity of this problem is unknown, even for $k=2$, and we view this as an interesting open problem. More precisely, suppose two items are hidden in $n$ boxes with search costs according to a known distribution, so that the probability they are in boxes $i$ and $j$ is some $p_{ij}$, for $i \le j$. What ordering of the boxes minimizes the total expected cost of finding both objects?


\subsection{Overview of Paper}
\label{sec:overview}

We take two approaches to using an $\alpha$-approximate best response oracle to find approximate solutions to games. In Section~\ref{sec:ellipsoid}, we present an algorithm 
using an existing {\em ellipsoid approach} from the approximation algorithms literature.
This approach uses the ellipsoid algorithm with an approximate separation oracle
to achieve approximately optimal solutions to an LP and to its dual.
The approach was previously used in work on solving generic LPs, and was used to develop approximation algorithms for a number of different combinatorial optimization problems, including routing, scheduling, and graph coloring
(\cite{Jansen03}, ~\cite{CV02}, ~\cite{FS12} and ~\cite{FKN12}). However, while the ellipsoid algorithm is theoretically important as the first polynomial time algorithm for solving LPs, it is extremely slow and is not used in practice. Thus we include our algorithm from Section~\ref{sec:ellipsoid} for theoretical interest only.

In Section~\ref{sec:mult-weights} we use a {\em multiplicative weights update} approach: a good general survey of this method can be found in \cite{AHK}. Approximately 20 years ago, ~\cite{FS99,FS96} gave a multiplicative update algorithm that computes
near-optimal strategies for two-player zero-sum games.   
Their algorithm yields strategies with payoffs that are within an {\em additive} $\varepsilon$ of 
the game value $V$.  Specifically, the strategy for Player I guarantees an expected payoff of at least $V-\varepsilon$, and the strategy for 
Player II guarantees an expected payoff of at most $V+ \varepsilon$.
The algorithm assumes that the payoffs are in the interval $[0,1]$, so
$\varepsilon$ represents a fraction of the range $[0,1]$ of payoffs.

Surprisingly, the approach of Freund and Schapire 
does not seem to have been extended to yield 
multiplicative approximations: that is, strategies that guarantee an expected payoff within a multiplicative factor $(1+\varepsilon)$ of the value. Multiplicative approximations are more common in the literature on approximation algorithms. It is also easy to show that for small $\varepsilon$, our notion of $(1+\varepsilon)$-optimal is equivalent to the notion of a {\em relative $\varepsilon$-approximate Nash equilibrium} from \cite{D13}, in the case of zero-sum games (see the Appendix).

In Section~\ref{sec:mult-weights}, we obtain $(1+\varepsilon)$-optimal strategies.  
We use a multiplicative update
rule of~\cite{GK07}, and we rely somewhat on their analysis, though their work does not automatically imply the existence of approximately optimal strategies for Player I. We therefore combine this with a variation of the analysis of~\cite{FS99}
as presented by~\cite{AHK}. We make some changes, 
due to the form of our multiplicative approximation factor,
and the fact that our algorithm is designed to accommodate an approximate
best response oracle. However, our primary contribution is a bringing together of ideas resulting in 
what we believe is the first self-contained exposition of an algorithm and its analysis for computing multiplicative approximations of optimal strategies in a zero-sum game.

We then apply our algorithms to a variety of games in Section~\ref{sec:app}.

In Section~\ref{sec:experiments} we describe the results of experiments we performed, in which we compare our multiplicative weights algorithm and the multiplicative weights algorithm of \cite{FS99} on the game BOX. 


\section{Preliminaries}
\label{sec:prelim}
In what follows, we consider a zero-sum game between a maximizing Player I with $n$ pure strategies and a minimizing Player II with a number of strategies that 
that could be as large as $2^{p(n)}$, for some polynomial $p$.
The payoff function is $C$ and all the payoffs are positive. We assume there 
is an $\alpha$-approximate best response oracle, for some constant $\alpha \geq 1$. That is, 
there is a polynomial time algorithm (polynomial in $n$) which, for any given mixed 
strategy $\mathbf{x}$ of Player I, computes a pure strategy $j_{\alpha}(\mathbf{x})$ 
such that $C(\mathbf{x},j_\alpha(\mathbf{x})) \le \alpha \min_i C(\mbfx,i)$. In other words, the oracle computes a pure strategy for Player II whose payoff against $\mathbf{x}$ is within a factor $\alpha$ of the payoff given by Player II's best response.

Let $V^*$ be the value of the game.
We say that a strategy $\hat{\mbfx}$ for Player I is
$\alpha$-optimal if
$C(\hat{\mbfx},\mbfy) \geq V^*/\alpha$ for any strategy $\mbfy$ of Player II.
We say a strategy for Player II is
$\alpha$-optimal if
$C(\mbfx,\hat{\mbfy}) \leq \alpha V^*$ for any strategy $\mbfx$ of Player I.


We assume that mixed strategies $\mbfx$ and $\mbfy$ for the two players
are represented in sparse form, as a list
containing the index of each non-zero entry, together with its value.
This is especially important for $\mbfy$, as the number of pure strategies
for Player II may not be polynomial in $n$, and we want the
running time to depend polynomially on $n$.

Let $A$ denote the game matrix, so that $A[i,j] = C(i,j)$.
We assume an oracle for computing $A[i,j]$.
That is, we assume there
is a polynomial-time algorithm (polynomial in $n$) that takes as input
(representations of) pure strategies $i$ and $j$
of Players I and II, and outputs $A[i,j]$.

\section{Ellipsoid Approach}
\label{sec:ellipsoid}

Let $\mu=\max_{i,j} A[i,j]$ be the largest entry in the game matrix.
We assume this quantity is given as part of the input,
or that it can be computed in time polynomial in $n$.
In fact, it is sufficient to compute an upper bound $\mu$ that is polynomial in
$\max_{i,j} A[i,j]$ and $n$. We assume in this section that the payoffs are integers, since the standard bounds on the runtime of the ellipsoid algorithm assume that the input LP has integer coefficients.

Consider the standard LP representing the problem of
finding an optimal strategy for Player I.

\vspace{6pt}
\hrule
\vspace{6pt}
{\bf \noindent LPI:}
Maximize $V$\\[3pt] 
such that\\[3pt]
\mbox{\ \ \ \ }(1) $\sum_i x_i A[i,j] \geq V \mbox{ for all } j$\\[3pt] 
\mbox{\ \ \ \ }(2) $\sum_i x_i = 1$\\[3pt]
\mbox{\ \ \ \ }(3) $x_i \geq 0 \mbox{ for all } i$\\[3pt] 
\vspace{6pt}
\hrule
\vspace{16pt}

A separation oracle for LPI takes as input an assignment
to the variables of the LP, and either reports that it is a feasible solution,
or returns a hyperplane separating 
the assignment from the feasible region.
If we had a separation oracle for LPI, we could use the ellipsoid
algorithm to obtain an optimal strategy for Player I.

Although we do not have such an oracle,
we can use the $\alpha$-approximate best response oracle
as an approximate separation oracle
for LPI.
A query to a separation oracle for LPI
corresponds to a pair $<\mbfx,V>$, where the first element is the assignment to the $x_i$
and the second is the assignment to $V$.
We can answer query $<\mbfx,V>$, approximately, using the following procedure:

\bigskip 

\begin{enumerate}
\item Check if a constraint of type (2) or (3) of LPI is violated by
assignment $<\mbfx,V>$.  If so, return such a violated constraint (separating hyperplane) as the answer to the query and exit. 
\item Otherwise, $\mbfx$ corresponds to a mixed strategy of Player I.  Query the
$\alpha$-approximate best response oracle on $\mbfx$.
Let $j^*=j_{\alpha}(\mathbf{x})$ be the 
returned strategy for Player II. 
\item Check whether $\sum_i x_i A(i,j^*) \geq V$ is satisfied by the queried assignment $<\mbfx,V>$.
\item If not, return the violated constraint $\sum_i x_i A(i,j^*) \geq V$ as the answer to the query.
Otherwise, answer the query by reporting (possibly incorrectly) that
$<\mbfx,V>$ is a feasible solution.
\end{enumerate}

\bigskip

Note that in the last step $<\mbfx,V>$ may be returned, even though it is not necessarily feasible.
However, in this case $<\mbfx,V/\alpha>$ is feasible,
because $j^*$ was returned by the $\alpha$-approximate best response oracle.


With this approximate separation oracle,
we can apply the ellipsoid approach
from the approximation algorithms literature,
referenced above in
Section~\ref{sec:overview}. (In fact,~\cite{FKN12} uses a 
minor variant of the approach that yields a somewhat smaller restricted dual LP. For simplicity of presentation, we do not discuss this variant here.)

It yields an algorithm that computes $\alpha$-optimal strategies for Players I and II.
The algorithm makes use of the dual of LPI, which computes an optimal strategy for Player II.
We present the dual LP as LPII.

\vspace{6pt}
\hrule
\vspace{6pt}
{\bf \noindent LPII:}
Minimize $V$\\ 
such that\\[3pt]
\mbox{\ \ \ \ }(1) $\sum_j y_j A[i,j] \leq V \mbox{ for all } i$\\[3pt] 
\mbox{\ \ \ \ }(2) $\sum_j y_j = 1$\\[3pt]
\mbox{\ \ \ \ }(3) $y_j \geq 0 \mbox{ for all } j$\\[3pt] 
\vspace{6pt}
\hrule
\vspace{16pt}

The algorithm is as follows.
To compute an $\alpha$-approximate strategy for Player I, run the ellipsoid algorithm
using the approximate separation oracle for LPI described above.
Let $S$ be the set of violated constraints returned by the oracle.
Let $\hat{\mbfx}$ denote the strategy that is returned.  
Return it as the $\alpha$-optimal strategy for Player I.

Let $S'$ be the set of variables $y_j$ in LPII, corresponding to
the constraints in $S$.
Let LPIIRestrict denote the restriction of LPII that is produced
by setting all $y_j \not\in S'$ to 0.
Generate LPIIRestrict explicitly, computing $A[i,j]$
for all $j \in S'$ and for all $i$, to obtain the constraints.
Solve this LP using a standard polynomial-time LP algorithm,
such as the ellipsoid algorithm, to obtain an optimal solution $\hat{\mbfy}$
for LPIIRestrict.
Return $\hat{\mbfy}$ as the strategy for Player II.

This algorithm is the basis for the following theorem.
In the proof, we provide the analysis of the algorithm for
the benefit of the reader; the analysis is essentially the same as that presented by~\cite{Jansen03} for solving general LPs and their duals.

\begin{theorem}\label{thm:ellipsoid} Let $V^*$ be the value of the game. There exists 
an algorithm that computes strategies $\hat{\mbfx}$ and $\hat{\mathbf{y}}$ 
	for Players I and II 
	such that
	\begin{align*}
	C(\mbfx,\hat{\mathbf{y}}) &\le \alpha  V^* \mbox{ for any strategy $\mathbf{x}$ of Player I} \mbox{ and} \\ C(\hat{\mbfx},{\mathbf{y}}) & \ge \frac {V^*} {\alpha} \mbox{ for any strategy $\mathbf{y}$ of Player II.}
	\end{align*}
This algorithm, which uses the ellipsoid algorithm,
	runs in time polynomial in $n$ and $\ln {\mu}$, where $\mu$ is the largest
entry in the game matrix.
\end{theorem}

\begin{proof}
The algorithm presented above computes the strategies
$\hat{\mbfx}$ 
and $\hat{\mbfy}$. 
We now prove that $\hat{\mbfx}$ 
and $\hat{\mbfy}$ are
$\alpha$-optimal strategies for Players I and II respectively.
To show that
$\hat{\mbfx}$ is an $\alpha$-optimal strategy for Player I,
let LPI($v$) denote the feasibility problem for
the set of constraints of LPI, with $V=v$.
Let $P(v)=\{\mbfx~|~\mbfx$ is a feasible solution to LPI($v$)$\}$.
First consider what happens if you run the
ellipsoid algorithm using an exact separation oracle for LPI.
It performs a binary search
on the values in an interval $I := \{0,\gamma, 2\gamma, \ldots, \mu\}$, 
where $\gamma$ and $1/\gamma$ is polymomial
in $n$ and $\ln \mu$.
The optimal value of LPI is guaranteed to be
a value in $I$.
During the binary search, for each tested value $v$ in $I$,
the ellipsoid algorithm makes a sequence $Q_v$ of queries $<\mbfx,V>$ to the separation oracle, where
in each query, $V=v$.  The querying continues until either
(a) the separation oracle returns the queried assignment $<\mbfx,v>$  (indicating that $\mbfx \in P(v)$) or
(b) the set of violated constraints returned by the separation oracle, in answer to the queries in $Q_v$,
implies that $P(v) = \emptyset$ (because there is no assignment satisfying
all constraints in that set).
The binary search 
finds the largest value of $v \in I$ for which $P(v) \neq \emptyset$.
The ellipsoid algorithm then returns this largest value $v$,
along with the associated $\mbfx \in$ LPI($v$), as the optimal solution to LPI.

Now consider what happens when the ellipsoid algorithm is
instead run with
the approximate separation oracle described above.
Recall that $V^*$ is the optimal value of LPI.
The binary search is run, and ends by returning a final value of $v$,
with an associated assignment $\mbfx$.
Let $v^f$ denote this final value, and let $\hat{\mbfx}$ be the associated assignment.
By the properties of the approximate separation oracle, 
and the behavior of binary search,
$\hat{\mbfx} \in P(v^f/\alpha)$ and $P(v^f+\gamma) = \emptyset$.
Because $P(v^f+\gamma) = \emptyset$,
\begin{equation}
\label{vstarvf}
v^f \geq V^*. 
\end{equation}

Further, because $\hat{\mbfx} \in P(v^f/\alpha)$, it follows that
for all strategies $\mbfy$ of Player II, 
\begin{equation}
\label{xbound1}
C(\hat{\mbfx},\mbfy) \geq v^f/\alpha,
\end{equation}
and hence from (\ref{vstarvf}),
\begin{equation}
\label{xbound}
C(\hat{\mbfx},\mbfy) \geq V^*/\alpha 
\end{equation}
for all strategies $\mbfy$ of Player II. Thus
$\hat{\mbfx}$ is an $\alpha$-optimal strategy for Player I.

By the well-known results in~\cite{GLS} on
the runtime of the ellipsoid algorithm, 
the above computation of $\hat{\mbfx}$ takes time polynomial in $n$ and $\ln \mu$.

The proof that $\hat{\mbfy}$ is
an $\alpha$-approximate strategy for Player II is as follows.
Let $V^R$ denote the optimal value of LPIIRestrict.
Since $\hat{\mbfy}$ is an optimal solution to LPIIRestrict, 
\begin{equation}
\label{vropt}
C(\mbfx, {\hat{\mbfy}}) \leq V^R
\end{equation}
for all strategies $\mbfx$ of Player I.

Because $\hat{\mbfx} \in P(v^f/\alpha)$,
\begin{equation}
\label{vf}
v^f/\alpha \leq V^*.
\end{equation}
Since LPII is dual to LPI,
the optimal value of LPII is also $V^*$.

Let LPIRelax denote the dual of LPIIRestrict obtained from LPI
by removing the constraints 
$\sum_i x_i A[i,j] \geq V$ for $j \not\in S'$.
By duality, $V^R$ is also the optimal value for LPIRelax.
The violated constraints returned when running the ellipsoid algorithm with the
approximate separation oracle for LPI are also violated constraints
of LPIRelax.  Therefore, the answers given by the approximate separation oracle
for LPI could also have been given for LPIRelax. 
Since the execution of the ellipsoid algorithm depends entirely on
the answers to the separation oracle queries, and not
on other properties of the LP, 
it follows by the same argument used to prove (\ref{vstarvf}) that
\begin{equation}
v^f \geq V^R. 
\end{equation}
Combining that with (\ref{vf}), we get
$V^R \leq \alpha V^*$.
From (\ref{vropt}), it  immediately follows that $\hat{\mbfy}$ is an
$\alpha$-optimal strategy for Player II.

Because the computation of $\hat{\mbfx}$ takes time polynomial in $n$ and $\ln{\mu}$,
the number of constraints in LPIIRestrict is polynomial in those parameters.
It follows that computing
$\hat{\mbfy}$ takes time polynomial in those quantities.
\end{proof}

\section{Multiplicative Weights Approach}
\label{sec:mult-weights}
We now present an alternative approach to solving (or approximately solving) games using an $\alpha$-approximate best response oracle. This approach is more ``combinatorial'' in nature, in the sense that it does not depend on the values of the payoffs in the same way as the approach of Section~\ref{sec:ellipsoid}.

Our algorithm takes the same high level approach as the algorithm of \cite{FS99}, which finds strategies that approximate the optimal strategies within an {\em additive} $\varepsilon$.
Player I begins with a uniformly random strategy. 
Player I then repeatedly plays mixed strategies, to which Player II replies with best responses. In each round, 
Player I updates her strategy, using some multiplicative update rule based on the best response of Player II. 
The more successful the pure strategy is (that is, the higher payoff it achieves), the more weight Player I gives it 
in the next round.   

The algorithm is based on a simple observation: finding optimal strategies for Player II is easily reducible to the problem of solving a {\em Packing LP}.
An approximately optimal solution to the latter problem is given by a multiplicative weights update algorithm, due to~\cite{GK07}.
Our algorithm uses essentially the same update rule, and we also rely on part of their analysis to show that this gives an approximately optimal solution for Player II. We show how the multiplicative weights update algorithm also yields an approximately optimal strategy for Player I. We believe the strategy for Player I to be an original contribution, which does not immediately follow from the analysis of \cite{GK07}.

The relation between our games and packing LPs is as follows.
A packing LP is a linear program of the form $\max\{\mathbf c^T\mathbf y~|~A\mathbf y \leq \mathbf b, \mathbf y \geq 0 \}$
where $A, \mathbf b$, and $\mathbf c$ all have positive entries.
If $A$ is the payoff matrix of a two-player zero-sum game, 
$\mathbf b^T = (1, \ldots, 1)$, $\mathbf c^T = (1, \ldots, 1)$, and $\mathbf{y^*}$
is the optimal solution to the LP, it is easy to show that 
the strategy that chooses each column $j$ with
probability proportional to $y^*_j$ is an optimal strategy for the Player II.
Thus one can reduce the problem of computing an optimal strategy for Player II
to the problem of solving the packing LP.\  
Because the reduction is so simple
(just a scaling of the variables), we integrate it into our algorithm, rather than
performing it explicitly as a separate step.
We note that a special case of this reduction was exploited previously by~\cite{CDHW} in the solution
of a variant of the game BOX. However, the resulting specialized
packing LP was then solved by different methods.

We first describe our algorithm in detail, then give a bound on the run time, and finally prove that it works. 

Fix $\eta > 0$ and $\delta >0$ (we will specify their values later) and let $\mbfx^{(1)}$ be the mixed strategy for Player I given by $x^{(1)}_i = 1/n$ for all $i$. 

\bigskip

\begin{enumerate}
	\item Set $t=1$.
	\item Use the $\alpha$-approximate best response oracle to compute $j^{(t)} = j_\alpha(\mbfx^{(t)})$.
	\item Set $M^{(t)} = \max_i C(i,j^{(t)})$. (Positive payoffs guarantee that $1/M^{(t)}$ is well defined.)
	\item Define $x^{(t+1)}_i = \left( \frac{1 + \eta C(i,j^{(t)})/M^{(t)}}{1 + \eta C(\mbfx^{(t)},j^{(t)})/M^{(t)} } \right) x^{(t)}_i$ for all $i$. (The expression in the denominator is a normalizing factor to ensure that $\mbfx^{(t+1)}$ is a mixed strategy.)
	\item If $f(t):=  \prod_{t' \le t} \left(1 + \eta C(\mbfx^{(t')},j^{(t')}) / M^{(t')} \right) > 1/(\delta n)$, then set $T=t$ and stop. Else increase $t$ by $1$ and return to Step 2.
\end{enumerate}

\bigskip

Let $\hat{\mbfy}$ be the mixed strategy that chooses pure strategy $j^{(t)}$ with probability proportional to $1/M^{(t)}$. More precisely, $\hat{y}_j$ is proportional to $\sum_{\{t:j^{(t)}=j\}} 1/M^{(t)}$. Let $\hat{t}$ be the value of $t$ that maximizes $C(\mbfx^{(t)}, j^{(t)})$ and let $\hat{\mbfx} = \mbfx^{(\hat{t})}$. We will later show that these strategies $\hat{\mbfx}$ and $\hat{\mbfy}$ are approximately optimal.

The proof of the following lemma is from~\cite{GK07} (following the presentation in \cite{LPRY08}).

\begin{lemma} \label{lem:Tbound}
	The number of iterations $T$ of the algorithm satisfies
	\begin{align}
	T \le n \left( 1 - \frac{\ln \delta}{\ln(1+\eta)} \right). \label{eq:Tbound}
	\end{align}
\end{lemma}

\begin{proof}
For $i=1,\ldots,n$, let $z^{(1)}_i = \delta$ and for $t \ge 2$ let $\mbfz^{(t)} = \delta n f(t-1) \mbfx^{(t)}$. It follows that 
\begin{align}
z^{(t+1)}_i = (1 + \eta C(i, j^{(t)})/M^{(t)}) z^{(t)}_i. \label{eq:zrecur}
\end{align}
Note that the stopping condition $f(t) > 1/(\delta n)$ is equivalent to the condition $\sum_i z^{(t+1)}_i > 1$, since
\[
\sum_i z^{(t+1)}_i = \sum_i \delta n f(t) x^{(t+1)}_i = \delta n f(t).
\]
Therefore, since the algorithm terminates on the $T$th iteration, we must have $\sum_i z^{(T)}_i \le 1$ and $\sum_i z^{(T+1)}_i > 1$. 

In each iteration of the algorithm, there is some pure strategy $i$ for Player I such that $C(i,j^{(t)}) = M^{(t)}$. Call such a strategy a {\em bottleneck strategy}. If $i$ is a bottleneck strategy in iteration $t$ then $z_i^{(t+1)} = (1+\eta) z_i^{(t)}$, so if $i$ is a bottleneck strategy a total of $m$ times then $z^{(T+1)}_i \ge \delta (1+\eta)^m$. But also $z^{(T+1)}_i \le 1+\eta$, since otherwise $z^{(T)}_i > 1$, contradicting $\sum_i z^{(T)}_i \le 1$. Putting together these inequalities,
\begin{align*}
\delta(1+\eta)^m &\le z^{(T+1)}_i \le 1+\eta, \mbox{ so }\\
m &\le 1 - \frac{\ln \delta}{\ln(1+\eta)}.
\end{align*}
Inequality~(\ref{eq:Tbound}) follows from the fact that Player I has $n$ pure strategies and in each iteration at least one of them is a bottleneck strategy.
\end{proof}

We now prove that for an appropriate choice of the parameters $\delta$ and $\eta$, the strategies produced by the algorithm are approximately optimal and that the total runtime is not too large.

\begin{theorem}\label{thm:main} Let $V^*$ be the value of the game. For any $\varepsilon >0$, there exist algorithms that compute strategies $\hat{\mbfx}$ and $\hat{\mathbf{y}}$ 
	for Player I and II 
	such that
	\begin{align*}
	C(\mbfx,\hat{\mathbf{y}}) &\le \alpha (1+\varepsilon) V^* \mbox{ for any strategy $\mathbf{x}$ of Player I} \mbox{ and} \\ 
	C(\hat{\mbfx},{\mathbf{y}}) & \ge \frac {V^*} {\alpha(1+\varepsilon)} \mbox{ for any strategy $\mathbf{y}$ of Player II.}
	\end{align*}
	The algorithm that computes $\hat{\mbfy}$ runs in time polynomial in $n$ and $1/\varepsilon$; the algorithm that computes $\hat{\mbfx}$ runs in time polynomial in $n$, $1/\varepsilon$ and $\alpha$.
\end{theorem}

\begin{proof}

For each $i$, we have
\[
x^{(T+1)}_i = \left( \prod_{t \le T} \frac{1 + \eta C(i,j^{(t)})/M^{(t)}}{1 + \eta C(\mathbf{x}^{(t)},j^{(t)})/M^{(t)}} \right) x^{(1)}_i 
=    \left(  \prod_{t \le T} \frac{1 + \eta C(i,j^{(t)})/M^{(t)}}{1 + \eta C(\mathbf{x}^{(t)},j^{(t)})/M^{(t)}} \right) \frac 1 n .
\]	
Since every coordinate of $\mbfx^{(T+1)}$ is at most $1$, we have
\begin{align}
n \ge  n x^{(T+1)}_i = \frac{\prod_{t \le T} 1 + \eta C(i,j^{(t)})/M^{(t)}}{\prod_{t \le T}1 + \eta C(\mathbf{x}^{(t)},j^{(t)})/M^{(t)}} 
\ge \frac{ (1+ \eta)^{\sum_{t \le T} C(i, j^{(t)})/M^{(t)}}}{\exp(\eta \sum_{t \le T} C(\mbfx^{(t)},j^{(t)})/M^{(t)}))} \label{ineq}
\end{align}
The second inequality above follows from the facts that
\[
(1+\eta z) \ge (1+\eta)^z \mbox{ for } z \in [0,1], 
\]
and
\[
(1+ z) \le \exp(z) \mbox{ for all }z.
\]
Taking natural logarithms in~(\ref{ineq}), we get
\[
\ln n \ge \ln(1+\eta) \sum_{t \le T} C(i, j^{(t)})/M^{(t)} - \eta \sum_{t \le T} C(\mbfx^{(t)},j^{(t)})/M^{(t)}.\\
\]
Now we use the fact that $\ln(1+\eta) \ge \eta - \eta^2$ 
and rearrange:
\[
\sum_{t \le T} C(\mbfx^{(t)}, j^{(t)})/M^{(t)} \ge  - \frac{\ln n}{\eta} + (1-\eta)\sum_{t \le T} C(i, j^{(t)})/M^{(t)}.
\]
By the linearity of the cost function and the right-hand side, the inequality above is still true if we replace the terms $C(i,j^{(t)})$ by $C(\mbfx,j^{(t)})$, where $\mbfx$ is any arbitrary strategy for Player I. Combining this with the fact that $C(\mbfx^{(t)}, j^{(t)}) \le \alpha V^*$, we get
\begin{align}
\alpha V^* \sum_{t \le T} 1/M^{(t)} \ge \sum_{t \le T} C(\mbfx^{(t)}, j^{(t)})/M^{(t)} \ge - \frac{\ln n}{\eta} + (1 -\eta) \sum_{t \le T} C(\mbfx, j^{(t)})/M^{(t)}. \label{ineq2}
\end{align}

By definition of $\hat{\mathbf{y}}$,
\[
C(\mbfx,\hat{\mathbf{y}}) = \frac {\sum_{t \le T} C(\mathbf{x},j^{(t)})/M^{(t)}}{\sum_{t \le T} 1/M^{(t)}},
\]
and substituting this into~(\ref{ineq2}) and rearranging gives
\begin{align}
\alpha V^*  \ge - \frac{\ln n}{\eta \sum_{t \le T} 1/M^{(t)}} + (1 -\eta) C(\mbfx,\mathbf{\hat{y}}). \label{ineq3}
\end{align}
By the stopping condition, $f(T) > 1/(\delta n)$, so
\begin{align}
\frac{1}{ \delta n} &<  \prod_{t \le T} \left(1 + \eta C(\mbfx^{(t)},j^{(t)}) / M^{(t)} \right) \nonumber \\
& \le  \prod_{t \le T} \left( 1 + \eta \alpha V^* / M^{(t)} \right) \nonumber \mbox{ (by definition of $j^{(t)}$)}\\
& \le \exp \left(\sum_{t \le T} \eta \alpha V^*/M^{(t)} \right), \mbox{ so taking logs and rearranging,} \nonumber\\
\sum_{t \le T} \frac{1}{M^{(t)}} & \ge -\frac{\ln(n \delta)}{\eta \alpha V^*}. \label{ineq5}
\end{align}
We will see later than $n \delta <1$, so that $\ln(n \delta)$ is negative. Combining~(\ref{ineq5}) with~(\ref{ineq3}) gives
\begin{align}
C(\mbfx,\hat{\mbfy}) \le   \frac{ \alpha V^* }{(1 - \eta)\left(1+\ln n /\ln \delta \right)}  . \label{eq:Cbound}
\end{align}
Let $\eta = \eta_0 = ((1+\varepsilon)^{1/2} - 1)/2$ and let $\delta = n^{-1/\eta}$, and use $(1-\eta)^{-1} \le (1+2\eta)$ to get
\[
C(\mbfx,\hat{\mbfy}) \le  \frac{ \alpha V^*}{(1 - \eta)^2} \le (1+2\eta)^2 \alpha V^* = \alpha(1+\varepsilon)V^*.
\]
With these choices of $\eta$ and $\delta$, by Lemma~\ref{lem:Tbound}, the number of iterations $T$ of the algorithm satisfies
\begin{align*}
T &\le n \left(1 + \frac{(1/\eta)\ln n}{\ln(1+\eta)} \right)\\
& \le n \left( 1+ (1/\eta+ 1/\eta^2)\ln n \right) \mbox{ (using $\ln(1+\eta) \ge \eta/(1+\eta)$)}\\
& \le n(1+(8/\varepsilon+ 8/\varepsilon^2) \ln n) \mbox{ (since $\eta \ge \varepsilon/8$)}.
\end{align*}
This is clearly polynomial in $n$ and $1/\varepsilon$. 

We now turn to Player I's strategy $\hat{\mbfx}$. We will choose different, smaller values for both $\eta$ and $\delta$ (to be specified later), giving a longer runtime for computing $\hat{\mbfx}$ that is polynomial in $n$, $1/\varepsilon$ {\em and} $\alpha$. If we wish to calculate both strategies we could either run the algorithm twice using the two different choices of $\eta$ and $\delta$, or simply run the algorithm once using the second choice. This would give Players II's strategy as well as Player I's since the right-hand side of~(\ref{eq:Cbound}) is increasing in $\delta$ and $\varepsilon$.

 First let $\mbfx^*$ be an optimal strategy of Player I. Then $C(\mbfx^*,j^{(t)}) \ge V^*$ because $\mbfx^*$ is optimal. Hence, substituting $\mbfx = \mbfx^*$ into~(\ref{ineq2}),
\begin{align}
\sum_{t \le T} C(\mbfx^{(t)}, j^{(t)})/M^{(t)} \ge - \frac{\ln n}{\eta} + (1 -\eta) \sum_{t \le T} C(\mbfx^*, j^{(t)})/M^{(t)} \ge - \frac{\ln n}{\eta} +  (1 -\eta) V^* \sum_{t \le T} 1/M^{(t)}. \label{ineq4} 
\end{align}
By definition of $j^{(t)}$, we have $ C(\hat{\mbfx}, j^{(\hat{t})}) \le \alpha C(\hat{\mbfx},\mathbf{y})$ for any strategy $\mathbf{y}$ so
\begin{align*}
C(\hat{\mbfx}, \mathbf{y}) &\ge \frac 1 \alpha C(\hat{\mbfx}, j^{(\hat{t})})\\
& \ge  \frac{1}{\alpha \sum_{t \le T} 1/M^{(t)}} \sum_{t \le T} C(\mbfx^{(t)},j^{(t)})/M^{(t)} \mbox{ (by definition of $\hat{t}$)}\\
& \ge \frac{-\ln n/\eta}{\alpha \sum_{t \le T} 1/M^{(t)}} + \frac{(1-\eta)V^*}{\alpha} \mbox{ (by~(\ref{ineq4}))} \\  
&\ge \frac{V^* \ln n}{\ln(n \delta)} + \frac{(1-\eta)V^*}{\alpha}  \mbox{ (by~(\ref{ineq5}))}\\
& = \frac{V^*}{\alpha} \left(1 - \eta + \frac{\alpha \ln n}{\ln(n \delta)}\right).
\end{align*}
Again, take  $\eta = \eta_0$ and this time take $\delta = \delta_0 = n^{-1-\alpha/\eta}$, so that
\[
C(\hat{\mbfx}, \mathbf{y}) \ge \frac{V^*}{\alpha}(1-2\eta)
\ge \frac{V^*}{\alpha(1+2 \eta)^2}
= \frac{V^*}{\alpha(1+\varepsilon)}.
\]
For these values of $\eta$ and $\delta$, the bound on the runtime given by Lemma~\ref{lem:Tbound} is
\begin{align*}
T &\le n \left(1 + \frac{(1+\alpha/\eta)\ln n}{\ln(1+\eta)} \right)\\
& \le n \left( 1+ (1/\eta+ 1/\eta^2)(\alpha+\eta)\ln n \right) \\
& \le n(1+(8/\varepsilon+ 8/\varepsilon^2)(\alpha+ \varepsilon) \ln n) \mbox{ (again using $\eta \ge \varepsilon/8$ and also $\eta \le \varepsilon$)}.
\end{align*}
This time $T$ is polynomial in $n$, $1/\varepsilon$ and $\alpha$. 
\end{proof}

We note that our algorithm uses $O(\frac{\alpha}{\varepsilon^2}(n\ln n))$ rounds of updates, while the algorithm of \cite{FS96} (which finds strategies within an {\em additive} $\varepsilon$ of optimal) has only a logarithmic dependence of $n$. In fact, a multiplicative approximation factor of $(1 + \varepsilon)$
cannot be achieved in $o(n)/\varepsilon^k$ rounds of a multiplicative update algorithm, for any $k$, as we state in the following proposition. The proof can be found in the appendix.

\begin{proposition} \label{prop:runtime}
	Let $k, \varepsilon >0$. There exist zero-sum games with $n$ strategies for each player, with positive payoffs, such that Player II has no $(1+\varepsilon)$-optimal strategy with support of size $o(n)/\varepsilon^k$. 
	\end{proposition}

\section{Applications to Search Games}
\label{sec:app}

In this section we describe how our algorithms can be applied to several known and new search games. Each game is played between a Hider, who corresponds to Player I in our theorems, and a Searcher, who corresponds to Player II. In this section we present several theorems concerning the existence of polynomial time algorithms for computing optimal or approximately optimal strategies for search games. These theorems refer to algorithms based on either the ellipsoid approach in Section~\ref{sec:ellipsoid} or the multiplicative weights approach in Section~\ref{sec:mult-weights}.

\subsection{Searching in Boxes with Precedence Constraints}
\label{sec:prec}
We begin with a generalization of the game BOX discussed in Subsection~\ref{sec:boxes}. In particular, we restrict the Searcher only to orderings that are consistent with some predefined precedence constraints on the boxes. That is, we suppose that there is a partial order $\prec$ on the boxes so that box $i$ can be searched before box $j \neq i$ if and only if $i \prec j$. The Hider's strategy set and the payoff function remain unchanged. We call this game PREC, and while it has not been studied before in the form we define it in here, a more general version of it was studied in \cite{FLV16}, as we will discuss in the next subsection.

In order to apply our algorithms to PREC we must consider the best response problem of choosing the search that minimizes the expected search cost for a given Hider distribution $\mbfx$. Similarly as for BOX, we can reframe this problem as the scheduling problem of finding a precedence-constrained ordering of a set of jobs with processing times and weights to minimize the sum of the weighted completion times. This problem is usually denoted $1|prec|\sum w_j C_j$ in the scheduling literature. It is well known to be NP-hard, and many 2-approximation algorithms can be found, for instance \cite{Chekuri} and \cite{Schultz}. These algorithms can be used as a 2-approximate best response oracle, so we obtain our first new result for search games.
\begin{theorem}
	There exist polynomial time algorithms that compute $2$-approximate strategies for both players in PREC.
\end{theorem}

\subsection{The Submodular Search Game}
\label{sec:submod}
This game was introduced by \cite{FKR16}, and further studied by \cite{FLV16}. As in the game, BOX, the Hider's strategy set is $[n]$ and the Searcher's strategy set is all permutations $\pi$ of $[n]$. The difference is in the payoff function. Let $f:S \rightarrow \mathbb{R}^+$ be a non-decreasing submodular set function on $S$, which corresponds to the cost of searching subsets of $S$. For a given hiding location $j$ and a given permutation $\pi$, the payoff $C(j,\pi)$ of the game is defined as the cost of the set of locations searched up to and including $j$. That is,
\[
C(j,\pi) = f(\{ i: i \le \pi^{-1}(j) \}).
\]
We call this game SUB. It is a further generalization of PREC. This can be seen by defining a cost function $f$ on subsets $S$ of $[n]$ in PREC such that $f(S)$ is the cost of all the elements of $[n]$ in the precedence-closure of $S$. Then, under this cost function, optimal strategies for SUB will correspond to optimal strategies in PREC. 

\cite{FLV16} also consider the best response problem for SUB, and they prove that this problem can be solved approximately, within a factor of $2$, generalizing the analogous result in the scheduling literature. It follows that our algorithms can be used to calculate $2$-approximate strategies for the players in SUB, and so we obtain our next new result for search games.
\begin{theorem}
	There exist polynomial time algorithms that calculate $2$-approximate strategies for both players in SUB.
\end{theorem}

\subsection{Expanding Search}
\label{sec:exp}

In PREC, the partial order $\prec$ on $[n]$ uniquely defines a directed acyclic graph, and if that graph is a tree then PREC is equivalent to the {\em expanding search game} on a tree, introduced in \cite{AL13b}. The game, which we will call EXP, can be played on any undirected graph with a distinguished vertex called the {\em root} and edges with costs corresponding to the time taken to traverse them. A pure strategy for the Hider is a vertex at which to hide, and a pure strategy for the Searcher is a sequence of edges, the first of which is incident to the root, and each other of which must be adjacent to a previously chosen edge, not necessarily the most recently chosen edge. This type of search on a graph is called an {\em expanding search}. The payoff is the search cost, which is the sum of the costs of the edges chosen by the Searcher up to and including the first edge that is incident to the Hider's chosen vertex. See \cite{AL13b} for motivations and applications of expanding search. 

For EXP played on a tree, the best response problem is what \cite{AL13b} call the {\em Bayesian problem} of finding a search that minimizes the expected payoff for a known hiding distribution on the vertices of the tree, and the authors give a solution to this problem. The solution to the equivalent scheduling problem $1|prec|\sum w_j C_j$ for tree-like precedence constraints is also known, and can be solved using the so-called {\em Sidney decomposition} proposed by \cite{Sidney}. Therefore, our algorithms can be applied to find optimal strategies for the players in the game. However, optimal strategies have already been found in closed form for this game in \cite{AL13b}.

If we consider EXP on a more general graph, then the best response problem can be shown to be NP-hard (\cite{Durr}). Whether the problem can be efficiently approximated is an open question. If this were the case then our results would imply that approximately optimal strategies for EXP could also be found.

\subsection{Expanding Search Ratio}
\label{sec:exp-r}
This is a game that was introduced by \cite{ADL16}. The game is played on a graph with $n$ vertices plus a root vertex and edges with costs. The Hider chooses one of the $n$ non-root vertices of the tree and the Searcher chooses an expanding search, as in the expanding search game considered in \cite{AL13b}. However, for a given vertex $v$ and a given expanding search, the payoff is not the search cost, but the ratio of the search cost to the shortest path distance $d(v)$ from the root to $v$. This payoff is called the {\em search ratio}.  It is a measure of the ``regret'' incurred by the Searcher in finding the Hider, and so we call the game EXPr.

\cite{ADL16} find a $5/4$-optimal strategy for the Searcher in EXPr for star graphs, which is generalized in \cite{ADL17} to general unweighted graphs (where all the edges have equal cost) and tree graphs. Independently, \cite{CDHW} have studied the same game for star graphs in the context of throughput maximization, and they give a full solution. But the solution of EXPr for any other class of graphs is unknown.

Consider EXPr played on a tree graph, and we will show that in this case, the best response problem can be solved in polynomial time. Indeed, suppose the Hider is located at vertex $v$ with probability $p(v)$. Let $x(v)=p(v)/d(v)$ and let $\hat{\mbfx}$ be the normalization of $\mbfx$, so that $\hat{\mbfx}$ is a probability vector. Then the best response problem is the problem of finding the expanding search that minimizes the expected search ratio with respect to $\mathbf{p}$, which is clearly equivalent to the problem of finding the expanding search that minimizes the expected {\em search time} with respect to $\hat{\mbfx}$. The latter problem is the best response problem for EXP, which, as discussed in Subsection~\ref{sec:exp}, can be solved in polynomial time. Therefore we obtain the following theorem.
\begin{theorem}
	There exist polynomial time algorithms that calculate optimal strategies for both players in EXPr played on a tree.
\end{theorem}

For other types of graphs, using our results depends on finding approximations to the best response problem for EXP.

\subsection{Search Games with Regret}
\label{sec:regret}
The payoff matrix of the game EXPr discussed in the previous subsection, can be obtained from the payoff matrix of the game EXP by dividing each row of the matrix by a constant, which is equal to the shortest path distance to the vertex corresponding to that row. If we take any search game, we can define a ``regret version'' of that game by dividing the row corresponding to each Hider pure strategy $j$ by some constant. There is usually a natural choice of constant, for example in SUB, we might divide the row corresponding to Hider strategy $j$ by $f(\{j\})$, the cost that would be incurred by the Searcher if she knew where the Hider was (provided $f(\{j\})>0$ for all $j$).

If there is an $\alpha$-approximate best response oracle for some search game, then it follows that there is also an $\alpha$-approximate best response oracle for any regret version of that game, using similar reasoning as in Subsection~\ref{sec:exp-r}. This implies, using our results, that we can find $\alpha$-optimal strategies in both games. So, for instance, we immediately find that for the regret version of SUB, there is an efficient algorithm that finds 2-approximate strategies. 

We sum this up below.

\begin{theorem}
	Suppose there is an $\alpha$-approximate best response oracle for some game we call GAME with payoff function $C$. Let GAMEr be the game defined by the payoff function $C'(i,j) := k_i C(i,j)$, where $k_1,\ldots,k_n$ are non-negative constants. Then there is an $\alpha$-approximate best response oracle for GAMEr, and hence there exist polynomial time algorithms for calculating $\alpha$-optimal strategies.
\end{theorem}

\subsection{Hide-Seek and Pursuit-Evasion}

We now discuss a different type of game, introduced by \cite{GalCasas} in the context of predator-prey interaction, though the game could equally apply to problems of national security such as the pursuit of a terrorist. As in the games discussed above, a Hider or hidden object is located in one of $n$ locations (this time without search costs). The Searcher can search a subset of locations size at most $k$, for some given $k$, and if the Hider's location $j$ lies in this set, a pursuit ensues and the Searcher captures him with probability $p_j$. The payoff of the game, which the Hider wishes to minimize and the Searcher to maximize, is the probability that the Hider is captured. \cite{GalCasas} provide a closed form solution to this game.

Here we introduce a generalization of this game, in which the locations have search costs $c_j$. A pure strategy for the Searcher is a subset of locations, the sum of whose search costs does not exceed $k$. The payoff, as before, is the probability the Hider is captured.  We call this game HSPE (hide-seek and pursuit-evasion). If all the locations' costs are equal to 1, then HSPE is equivalent to the original game of \cite{GalCasas}. Consider the best response problem, where the probability $x_j$ that the Hider is in location $j,j=1,\ldots,n$ is known and the problem is to choose a subset of locations that maximizes the probability of capture. Let $w_j = p_j x_j$. Then the problem is to choose a subset $S$ of $[n]$ of total cost $c(S):= \sum_{j \in S} c_j$ at most $k$, that maximizes $w(S):= \sum_{j \in S} w_j$. This is the classic Knapsack problem, which has a fully polynomial time approximation scheme (see \cite{Vazirani}). That is, given $\varepsilon>0$, there exists an algorithm that is polynomial in $n$ and $1/\varepsilon$ that approximates the solution of the knapsack problem within a multiplicative factor of $1+\varepsilon$. This implies the following. 
\begin{theorem}
	The problem of finding optimal strategies in the game HSPE has a fully polynomial time approximation scheme.
\end{theorem}

\subsection{An Infinite Game}
\label{sec:infinite}

We finish this section by discussing an infinite game studied by \cite{LS15}, for which our results do not directly apply because the Searcher has a strategy set of (countably infinite) cardinality. We are optimistic that our results could be extended to such games in future work, which is why we mention it here.

In the game of \cite{LS15}, a Hider is located in one of $n$ locations with search costs $c_j$ and capture probabilities $p_j$, as in the game, HSPE of the previous subsection. The Searcher chooses an infinite sequence of locations, and each time the Searcher examines the Hider's location $j$, she finds him independently with probability $p_j$. The payoff of the game is the Searcher's expected cost of finding the Hider for the first time.

As \cite{LS15} point out, the solution of the best response problem for this game is well known. Suppose the Hider is hidden in location $j$ with probability $x_j$. Then, as showed by Blackwell (reported in \cite{Matula}), the optimal policy for the Searcher is, at any time, to choose a location that maximizes the index $x_j p_j/c_j$, with the hiding probability $x_j$ being updated according to Bayes' Law after each search. Of course, because the optimal policy is an infinite sequence, it cannot necessary be concisely expressed, but an approximately optimal solution may be found, and \cite{LS15} exploit this fact to give an algorithm which they believe converges to an optimal search strategy in the game. We leave it for future work to develop a provably fully polynomial approximation scheme for the problems of determining optimal strategies in the game.

\section{Numerical experiments}
\label{sec:experiments}

In this section we test our multiplicative weights algorithm of Section~\ref{sec:mult-weights} and compare it to the algorithm of \cite{FS99}. The latter algorithm has the same form as ours, but updates the weights in each iteration using a different multiplicative factor.  It outputs strategies that ensure an expected payoff of at least $V^* - \varepsilon$ for Player I and at most $V^* + \varepsilon$ for Player II. The number of iterations is $O(\frac{\ln n}{\varepsilon^2})$, but it assumes the payoffs have been scaled to lie in $[0,1]$.

In order to compare the Freund and Schapire algorithm to ours, we note that theirs could be applied directly to obtain a multiplicative
approximation, but the resulting number of rounds would depend on the payoffs of the game.
More particularly, consider using their algorithm
to find a strategy for Player II that is within 
a {\em multiplicative} factor $1+\varepsilon$ of optimal.  
This requires scaling the payoffs to be in
$[0,1]$, which can be done by dividing them by $\mu$, the maximum entry
in the payoff matrix.  Achieving a multiplicative approximation of $1 + \varepsilon$ is 
equivalent to achieving an {\em additive} approximation of $\varepsilon V^*/\mu$, relative to the
scaled values. For Player I, achieving a optimal strategy that is within a multiplicative factor of $1+\varepsilon$ of optimal
means guaranteeing a payoff of at least $1/(1+\varepsilon)$ times the value of the game.
Thus achieving this approximation for Player I is equivalent to achieving an additive payoff of at least $1-\varepsilon'$ relative to the scaled values, where $\varepsilon' = \varepsilon/(1+\varepsilon)\cdot V^*/\mu$. However, we do not know the value of $V^*$, so in order to calculate the appropriate additive approximation, we have to use a lower bound $m$ for $V^*$. Thus, a multiplicative approximation of $1+\varepsilon$ is best compared to an additive approximate of $\varepsilon'$, where
\begin{align}
\varepsilon'=\frac{\varepsilon}{1+\varepsilon} \left(\frac{m}{\mu}\right). \label{eq:epsilon'}
\end{align}
So to ensure the same quality of approximation for {\em both} players we should use~(\ref{eq:epsilon'}) to calculate the appropriate additive approximation in the Freund and Schapire algorithm. The simplest choice for the lower bound $m$ would be the minimum value of any payoff in the game. This yields $O((\frac{\mu}{m \varepsilon})^2 (\ln n))$ iterations, a quantity that depends on the payoffs. If $\mu/m$ is constant, this algorithm requires fewer rounds than ours, asymptotically, but if $\mu/m$ is large compared to $n$ it could require an arbitrarily large number of rounds.

Here we test both our multiplicative weights algorithm and that of Freund and Schapire (as presented in \cite{AHK}) on instances of the game BOX. We randomly generated 40 different instances of the game, which we grouped into four sets. For Set 1 and Set 2, we took $n=5$, and for Set 3 and Set 4 we took $n=10$. The payoffs in Sets 1 and 3 were uniformly chosen at random to be an integer between $1$ and $10$, and in Sets 2 and 4 they were chosen uniformly at random between $1$ and $100$. The 40 instances of the game BOX can be found in Table~2 in the appendix.

We also considered four different values of $\varepsilon$. In order to fairly compare the two algorithms, for each value of $\varepsilon$ we used for the multiplicative approximate of $1+\varepsilon$ in our algorithm, we used $\varepsilon'$ for the Freund and Schapire algorithm, as given in~(\ref{eq:epsilon'}). For the lower bound $m$ on the value, we thought it unfair to take unthinkingly the smallest payoff in the game, and instead took $m$ to be the largest cost $\max_i c_i$ of the $n$ boxes. This is a better simple lower bound on the value, which the Hider can achieve by hiding the object in the highest cost box. 

An even better lower bound $m$ could be obtained by using the more sophisticated approach of taking the best response payoff $C(\mbfx, j(\mathbf x))$, where $\mathbf x$ is the uniform strategy $x_i = 1/n$ for all $i$ and $j(\mbfx)$ is Player II's best response to $\mbfx$. Therefore we then ran the Freund and Schapire algorithm a second time using this choice of $m$ (if it was larger).

The algorithms were implemented in Python and the experiments were run on a MacBook Pro laptop with a 3.1GHz Intel Core i7 processor and 16GB of RAM. The results are displayed in Table~\ref{tab:results}. Our algorithm is labeled ``HLP'', the Freund and Schapire algorithm with $m = \max_i c_i$ is labeled ``FS'' and the Freund and Schapire algorithm with the more sophisticated choice of $m$ as described above is labeled ``FS+''. We ran our algorithm, taking $\delta =\delta_0$ and $\eta=\eta_0$ (as defined in the proof of Theorem~\ref{thm:main}).

The number of iterations of our algorithm is determined by our theoretical results. That is, our theoretical results say that if you run the algorithm for a certain number of iterations, you are guaranteed to achieve a certain accuracy. But, in fact, you may be able to achieve that accuracy in fewer iterations in practice. So in our experiments, we explore both the number of iterations needed for the algorithm to terminate, and the number of iterations needed to achieve the desired accuracy. 

The columns in Table~\ref{tab:results} labeled ``Convergence time'' display what can be thought of as the average time the algorithms took to converge ``in practice''. More precisely, convergence time is the number of iterations until we could be sure that both strategies converged to the desired accuracy, so that the ratio of the maximum payoff guaranteed by Player I to the minimum payoff guaranteed by Player II is no more than $1+\varepsilon$. 

The columns in Table~\ref{tab:results} labeled ``Total time'' indicate the average number of iterations until the algorithms terminated. The time to produce each entry in the table (requiring running each algorithm on 5 or 10 different versions of the game) ranged from under one second to over 25 hours, depending on the setting of the parameters and the resulting number of iterations. We did not optimize our implementation for speed. To avoid excessive processing time in running our experiments, for $\varepsilon = 5\%$ and $\varepsilon = 10\%$, we did not attempt to run an algorithm to termination. If an algorithm exceeded 3 hours of processing time, we ensured it had  converged to the desired accuracy, and then terminated it prematurely. If any of the runtimes of the 10 instances of a set fell into this category, the corresponding entries in the table  are marked with an asterisk. 


In order to calculate the ``PI error'' and ``PII error'', for each experiment, we first calculated the value $V$ of each game, using the known closed form solution (see Subsection~\ref{sec:boxes}). Then, for Player I, the error is the smallest $p$ for which the final strategy output by the algorithm guarantees an expected payoff of at least $V/(1+p)$; for Player II the error is the smallest $p$ for which the final strategy guarantees an expected payoff of at most $V(1+p)$. 

The convergence times and total times are written to the nearest integer and the errors are written in percentages to two decimal places.

\begin{table}[h]
	\centering
	\caption{Table of results of numerical experiments.}
	\begin{tabular}{|c|c|ccc|ccc|ccc|ccc|}
		\hline
		&&\multicolumn{3}{|c|}{\textbf{Total time}} & \multicolumn{3}{|c|}{\textbf{Convergence time}} & \multicolumn{3}{|c|}{\textbf{PI error ($\%$)}} & \multicolumn{3}{|c|}{\textbf{PII error  ($\%$)}}\\
		$\varepsilon$ ($\%$) & Set & HLP & FS & FS+ & HLP & FS & FS+& HLP & FS & FS+& HLP & FS & FS+ \\
		\hline 
		\multirow{4}{*}{5} &1 & $16933$	& $28553$	& $11723$	& 970	&1394&	1033&	0.01&	0.00&	0.01	&0.28	&0.33	&0.42\\
		&2 & $16422$	& $22240$	& $13761$	&911&	1209&	1045&	0.01&	0.00&	0.00&	0.26&	0.33&	0.37\\
		&3& $*$	& $*$	& $*$&	1180&	3470&	1495 	&0.03&	0.01&	0.02&	0.60	&1.74&	0.78\\
		&4&$*$	& $*$	& $*$&	1176&	3421&	1496&	0.03&	0.01	&0.02&	0.61&	1.78&	0.78\\
		\hline
		\multirow{4}{*}{10} &1	&4359&	7835&	3217	&248&	367 &	272&	0.02&	0.01&	0.02&	0.55&	0.63&	0.81\\
		&2&	4232&	6103 	&3777&	233 &	320 &	276 &	0.02&	0.01&	0.02&	0.53&	0.62&	0.70 \\
		&3&	6943&	$*$&	6591&	302	&911 	&393&	0.06&	0.02&	0.04&	0.43&	0.91&	0.58 \\
		&4&	$6944 $&	$*$&	$6606$	&300&	897	&394 	&0.05&	0.02&	0.04&	0.42&	0.92&	0.57\\
		\hline
		\multirow{4}{*}{50} &1&	214 &	584& 	241 &	12	&21	&16&	0.15	&0.07&	0.11&	2.65&	2.37&	3.07\\
		&2&	209&	 455&	282 &	13	&21	&19 &	0.14&	0.06&	0.10&	2.40&	2.33&	2.64 \\
		&3&	340 &	2900&	492 &	15&	53&	24&	0.38&	0.07&	0.24&	2.12&	0.93&	2.24 \\
		&4&	340&	2894&	493 &	14&	53&	24&	0.36&	0.07&	0.24&	2.24&	0.94&	2.23\\	
		\hline
		\multirow{4}{*}{100} &1	&66&	261 &	108& 	2	&3&	3&	0.38&	0.11&	0.20	&4.76&	3.66&	4.72\\
		&2&	65&	203&	126	&3&	5&	5	&0.40	&0.11&	0.15&	4.53&	3.54&	4.05\\
		&3	&104	 &1290	&219&	3&	13	&7&	0.77&	0.12	&0.39&	4.20&	1.41&	3.45\\
		&4	&104 	&1287&	220 &	4&	15&	7&	0.75&	0.12&	0.38&	4.15&	1.43&	3.43\\
\hline
		\end{tabular}
				\label{tab:results}
\end{table}

We can see from Table~\ref{tab:results} that HLP ran much quicker on average than FS (in the cases that they both terminated) and ran quicker than FS+ for $\varepsilon=100\%$ and $50\%$, but not for $\varepsilon=10\%$ or $\varepsilon = 5\%$. In practice, HLP converged significantly more quickly on average even than FS+ for all four values of~$\varepsilon$.

Interestingly, by the time the algorithms stopped running, the actual error of the strategies was far lower than guaranteed by the theory, particularly for large $\varepsilon$. Indeed, for $\varepsilon=100\%$, for example, the average Player 1 error after the algorithms terminated was less than $1\%$ for each of the three algorithms, and for Player 2, was less than $5\%$. Comparing the total run times for $\varepsilon=100\%$ with the convergence times for $\varepsilon=5\%$, we see that in order to ensure an approximation ratio of $\varepsilon=5\%$, it was generally quicker to run the algorithm for $\varepsilon=100\%$. 



Comparing the average errors across the three algorithms, we note that the average Player I error for HLP was always at least as great as that of both FS and FS+. An obvious explanation is that the former usually ran for fewer iterations. In separate experiments, we ran HLP for a number of iterations equal to the maximum number of iterations performed by FS and FS+; in this case we found that the Player I error for HLP was comparable to, and usually smaller than that of FS and FS+.

For Player II, HLP had the smallest average error for both $\varepsilon = 5\%$ and $\varepsilon = 10\%$. For the larger values of $\varepsilon$, the average Player II error for HLP was comparable to that of FS+, and FS had the smallest average error, which again can be attributed to longer runtimes.

\section{Conclusion}
\label{sec:conclusion}

We have shown how we may use an algorithmic approach to find solutions or approximate solutions to search games, by exploiting oracles that find the best responses or approximate best responses. This has applications to the search games we have mentioned in Section~\ref{sec:app}, but we believe it may have further applications in other search games, and indeed in games studied in other fields of operations research such as security games.

\section*{Acknowledgements}
L. Hellerstein was partially supported by NSF Award IIS-1217968.

\section*{Appendix}

\subsection*{Equivalence of $(1+\varepsilon)$-Optimal Strategies and Relative $\varepsilon$-approximate Nash Equilibria for Zero-Sum Games}
We show here that for small enough $\varepsilon$, our notion of $(1+\varepsilon)$-optimal strategies is in a natural sense equivalent to the notion of relative $\varepsilon$-approximate Nash equilibria for zero-sum games with positive payoffs, from 
Daskalakis (2013).

Consider a (possibly non-zero-sum) two-player game, and for mixed strategies $\mathbf x$ and $\mathbf y$ of Players I and II, denote the payoffs to the two players by $C_1(\mathbf x, \mathbf y)$ and $C_2(\mathbf x, \mathbf y)$, respectively. A relative $\varepsilon$-approximate Nash equilibrium is a pair of strategies $(\mathbf{\hat{x}},\mathbf{\hat{y}})$ such that
\[
C_1(\mathbf{\hat{x}},\mathbf{\hat{y}}) \ge (1-\varepsilon) C_1(\mathbf x,\mathbf{\hat{y}}) \text{ and } C_2(\mathbf{\hat{x}},\mathbf{\hat{y}}) \ge (1-\varepsilon)C_2(\mathbf{\hat{x}},\mathbf y)
\]
for all strategies $\mathbf x$ of Player I and $\mathbf y$ of Player II.

Now suppose the game is zero-sum, so that $C:=C_1=-C_2$, and suppose some strategies $(\mathbf{\hat x},\mathbf {\hat y})$ define a relative $\varepsilon$-approximate Nash equilibrium. Let $V$ be the value of the game and let $(\mathbf{x^*},\mathbf{y^*})$ be any optimal strategies. Then for any Player II strategy $\mathbf{y}$,
\begin{align}
C(\mathbf{\hat{x}},\mathbf{y}) & \ge  \frac{1}{1+\varepsilon}C(\mathbf{\hat x},\mathbf{\hat y}) \label{eq1}\\
& \ge \frac{1-\varepsilon}{1+\varepsilon} C(\mathbf{x^*},\mathbf{\hat y}) \label{eq2}\\
& \ge \frac{1-\varepsilon}{1+\varepsilon} C(\mathbf{x^*},\mathbf{y^*}) \label{eq3}\\
& = \frac{1-\varepsilon}{1+\varepsilon} V \label{eq4} \\
& \ge \frac{V}{1+3 \varepsilon}, \nonumber
\end{align}
for all $\varepsilon >0$. Inequalities~(\ref{eq1}) and (\ref{eq2}) follow from the definition of a $\varepsilon$-approximate Nash equilibrium, inequality~(\ref{eq3}) follows from the optimality of $\mathbf{y^*}$ and equation~(\ref{eq4}) follows from the definition of the value. 

It can be similarly shown that $C(\mathbf{x},\mathbf{\hat{y}})  \le (1+ 3\varepsilon)V$ for any Player I strategy $\mathbf{x}$. So the strategies $(\mathbf{\hat x},\mathbf {\hat y})$ are $(1+3\varepsilon)$-optimal.

Now suppose the strategies $(\mathbf{\hat x},\mathbf {\hat y})$ are $(1+\varepsilon)$-optimal. Then for any Player I strategy $\mathbf x$,
\begin{align}
C(\mathbf{\hat x},\mathbf{\hat y}) & \ge \frac{V}{1+\varepsilon} \label{eq5}\\
& \ge \frac{1}{(1+\varepsilon)^2} C(\mathbf{x},\mathbf {\hat y}) \label{eq6}\\
& \ge (1- 3\varepsilon) C(\mathbf{x},\mathbf {\hat y}), \nonumber
\end{align}
for $\varepsilon \le 1/3$. Inequality~(\ref{eq5}) follows from the $(1+\varepsilon)$-optimality of $\mathbf{\hat x}$ and inequality~(\ref{eq6}) follows from the $(1+\varepsilon)$-optimality of $\mathbf{\hat y}$. 

Similarly, it can be shown that $C(\mathbf{\hat x}, \mathbf{\hat y}) \le (1+3 \varepsilon) C(\mathbf{\hat x},\mathbf y)$ for any Player II strategy $\mathbf y$. So the strategies $(\mathbf{\hat x},\mathbf {\hat y})$ form a relative $3 \varepsilon$-approximate Nash equilibrium for small enough $\varepsilon$.

\subsection*{Proof of Proposition 
	1}

\begin{proof}
	Consider the game whose $n \times n$ payoff matrix has the value $n$ in the diagonal entries, and the value $1$ in all the off-diagonal entries.
	Then the value 
	$V$ of the game is $2-1/n \leq 2$. Suppose Player II has a $(1+\varepsilon)$-optimal strategy $\mathbf y$ with support $o(n)/\varepsilon^k$. 
	Then asymptotically, Player II can ensure the payoff does not exceed $\max_{\mathbf x} C(\mathbf x, \mathbf y) \ge n \varepsilon^k/o(n)$.
	Since $\mathbf y$ is $(1+\varepsilon)$-optimal, it follows that
	\begin{align*}
	\frac{n\varepsilon^k}{o(n)} \le \max_{\mathbf x} C(\mathbf x, \mathbf y) \le (1+\varepsilon)V \leq (1+\varepsilon) \cdot 2. \label{eq:imposs}
	\end{align*}
	But this implies $\frac{\varepsilon^k}{2(1+\varepsilon)} \le \frac{o(n)}{n}$, a contradiction.
\end{proof}


\newpage

\subsection*{Data for Numerical Experiments}

\begin{table}[h]
	\centering
	\caption{Randomly generated instances of BOX game. Entries in the table correspond to costs.}
	\begin{tabular}{|c|c|c|c|c|c|c|c|c|c|c|}
		\hline
		& \multicolumn{10}{|c|}{Instances} \\
		&1&2&3&4&5&6&7&8&9&10\\
		\hline
		\multirow{5}{*}{Set 1} & 6&	7&	1&	3&	6&	6&	7&	8&	3	&1\\
		
		&1&	10&	2&	3&	7&	2	&8	&10	&8&	1\\
		
		&7&	4&	4&	10&	6&	7&	2&	9&	2&	1\\
		
		&4&	5&	4	&1	&2	&6	&9	&6&	1&	6\\
		
		&10	&6&	3	&6&	6&	6&	3	&8&	2&	3\\
		\hline
		\multirow{5}{*}{Set 2}& 57&	47&	44&	99	&22&	89&	95&	12	&10	&8\\
		
		&7&	42&	100	&10&	31&	66&	8	&52	&96	&32\\
		
		&81&	70	&67	&9&	79	&90	&12	&30	&1	&13\\
		
		&82	&4	&11&	32&	79	&17	&63	&62	&9	&32\\
		
		&48&	71&	5	&83	&45	&92	&32	&2	&59	&63\\
		\hline
		\multirow{10}{*}{Set 3} & 7&	2&	10&	6&	3&	8&	7&	8&	3&	2\\
		
		&5&	4&	2&	10&	1&	4&	7&	8&	5&	8\\
		
		&4	&4	&5	&9	&9	&3	&10&	7	&10&	5\\
		
		&5	&7	&2	&9	&6	&10	&5	&3	&6	&5\\
		
		&8	&2	&9	&3	&10&	9&	9&	5&	2&	8\\
		
		&6	&6	&6	&1	&4	&5	&8	&3	&1	&10\\
		
		&3	&1	&6	&7	&7	&7	&3	&3	&1	&4\\
		
		&4	&7	&6	&7	&6	&5	&8	&2	&10&	9\\
		
		&5	&2	&9	&10	&1	&5	&1	&2	&5	&1\\
		
		&6	&4	&1	&10&	5&	3&	7&	9&	1&	9\\
		\hline
		\multirow{10}{*}{Set 4} & 58	&51&	56&	11&	70&	15&	96&	54&	65&	62\\
		
		&12&	78&	99&	38&	83&	65&	71&	58&	41&	30\\
		
		&35	&100	&87	&51	&59	&99	&69	&65	&39	&96\\
		
		&35	&25	&87&	14&	84&	14&	84&	14&	7	&67\\
		
		&46	&47&	90&	8	&16	&98	&39&	98&	60&	60\\
		
		&2	&33	&45	&27&	15&	40&	61&	18&	90&	36\\
		
		&72	&96	&76	&57	&73	&90	&6	&15	&43	&25\\
		
		&91	&51	&66	&22	&60	&86	&91	&61	&77	&74\\
		
		&6	&28&	65&	27&	45&	53&	38&	57&	8	&98\\
		
		&54	&62	&17	&26	&88	&84	&86	&74	&5	&73\\
		\hline
	\end{tabular}
	\label{tab:data}
\end{table}

\end{document}